\documentclass[11pt,leqno]{article}
\usepackage{amssymb}
\usepackage{amscd}
\usepackage{amsmath}
\usepackage{amsthm}

\makeatletter
\def\l@section{\@tocline{1}{4pt}{1pc}{}{}}
\def\l@subsection{\@tocline{2}{0pt}{2pc}{5pc}{}}
\makeatother
\title{Decomposition and Parity of Galois Representations Attached to GL$(4)$}
\author{Dinakar Ramakrishnan\footnote{Partially supported by the NSF grant DMS-1001916;
dinakar@caltech.edu}}
\date{253-37 Caltech, Pasadena, CA 91125}
\begin{document}
\maketitle

\markboth{Dinakar Ramakrishnan}{Decomposition and Parity for GL$(4)$}
\markright{Decomposition and Parity for GL$(4)$}


\section*{Introduction}

Let $F$ be a number field, and $\pi$ an isobaric (\cite{dr:La}), algebraic (\cite{dr:Cl1}) automorphic representation of GL$_n({\mathbb A}_F)$. We will call $\pi$ {\it quasi-regular} iff at every archimedean place $v$ of~$F$, the associated $n$-dimensional representation $\sigma_v$ of the Weil group $W_{F_v}$ (defined by the archimedean local correspondence) is {\it multiplicity free}. For example, when $n=2$ and $F={\mathbb Q}$, a cuspidal $\pi$ is quasi-regular exactly when it is generated by a holomorphic newform $f$ of weight $\geq 1$. Recall that $\pi$ is {\it regular} iff at each archimedean $v$ the {\it restriction} of $\sigma_v$ to ${\mathbb C}^\ast$ is multiplicity free; hence any regular $\pi$ is quasi-regular, but not conversely.

When $F$ is totally real, an algebraic automorphic representation $\pi$ of GL$_{n}({\mathbb A}_F)$ is said to be {\it totally odd} iff it is odd at each archimedean place $v$, i.e., iff the difference in the multiplicities of $1$ and $-1$ as eigenvalues of complex conjugation in $\sigma_v$ is at most $1$; in particular if $n$ is even, these two eigenvalues occur with the same multiplicity. One readily sees that any quasi-regular $\pi$ is totally odd, but not conversely.
When $n=2$, $\pi$ is said to be {\it even} (or that it has even parity) at an archimedean place $v$ (of $F$) iff if it is not odd at $v$.

Still with $F$ totally real, let $c$ be one of the $[F:{\mathbb Q}]$ complex conjugations in the absolute Galois group ${{\mathfrak G}}_F=$Gal$(\overline F/F)$. Recall that an $n$-dimensional $\overline{\mathbb Q}_p$-representation $\rho$ of ${\mathfrak G}_F$ is {\it odd relative to $c$} if the trace of $\rho(c)$ lies in $\{1,0,-1\}$. It is {\it odd} if it is so relative to every $c$. When $n=2$, $\rho$ is said to be {\it even relative to $c$} if it is not odd relative to $c$, which is the same as the determinant of $\rho(c)$ being $1$.

One says that an isobaric, algebraic $\pi$ on GL$(n)/F$ has a fixed {\it archi\-medean weight} iff there is an integer $w$ such that at every archimedean place $v$ of $F$,
the restriction to ${\mathbb C}^\ast$ of $\sigma_v\otimes\vert\cdot\vert^{(n-1)/2}$ is a direct sum of characters $z \to z^{p_j}\overline{z}^{q_j}$ with $p_j+q_j =w$ ($\forall \, j$). Every {\it cuspidal} algebraic $\pi$ has such a weight (\cite{dr:Cl1}).
The first object here is to provide a proof of the following assertion, which was established for the regular (algebraic), cuspidal case in an earlier preprint \cite{dr:Ra1} which has remained unpublished:

\medskip

\noindent{\bf Theorem A} \, \it Let $F$ be a totally real number field, $n \leq 4$, $p$ a prime, $\pi$ a quasi-regular algebraic, isobaric automorphic representation of GL$_n({\mathbb A}_F)$ of a fixed archimedean weight, and $\rho$ an associated $n$-dimensional, Hodge-Tate $\overline{{\mathbb Q}_p}$-representation of ${\mathfrak G}_F$ whose local $L$-factors agree with those of $\pi$ (up to a shift) at almost all primes $P$ of $F$. Then the semisimplification $\rho^{\rm ss}$ (of $\rho$) does not contain any irreducible $2$-dimensional Galois representation which is even relative to some complex conjugation $c$.
\rm

\medskip

Since $\pi$ is {\it algebraic}, it is expected by a conjecture of Clozel (\cite{dr:Cl1} that there is an associated Galois representation $\rho$ (see also \cite{dr:BuzG}). When~$\pi$ is {\it regular and selfdual} and $F$ totally real, the existence of $\rho$ has been well known for some time for GL$(n)$ (\cite{dr:Cl2}, \cite{dr:PL}, \cite{dr:ClHLN}). Now it is also known for $\pi$ {\it non-selfdual} and regular, by the important recent work of Harris, Lan, Taylor and Thorne (\cite{dr:HLTT}). When $\pi$ is quasi-regular, the algebraicity is already in \cite{dr:BHR} (for $F$ totally real), where one finds the equivalent notion of a {\it semi-regular, motivically odd} cusp form, and much progress has been made in the construction of $\rho$ in the selfdual case in the thesis of Goldring (\cite{dr:Go1,dr:Go2}), when the base change of $\pi$ to a suitable $CM$ extension $K$ with $K^+=F$ is known to descend to a holomorphic form on a suitable unitary group associated to $K/F$.

Here are some remarks on the proof of Theorem A. When $\rho^{\rm ss}$ is a direct sum $\eta_{1}\oplus \eta_{2}$ with each $\eta_{j}$ two-dimensional, one sees easily that the summands
must have the same parity. But then we are left with the subtle task of ruling out
both of them being {\it even}. We first appeal to the exterior square construction of Kim (\cite{dr:K}), which supplies (using the Langlands-Shahidi method) an isobaric automorphic form $\Pi=\Lambda^2(\pi)$ on GL$(6)/F$, and then, more importantly, we make use of the fact that one has some control, thanks to the works of Shahidi (\cite{dr:Sh1,dr:Sh2}) and Ginzburg-Rallis (\cite{dr:GR}) concerning the analytic properties of the exterior cube $L$-function of $\Pi$, which has degree $20$. Luckily, this $L$-function is also related to the square of a twist of the symmetric square $L$-function of $\pi$, and we exploit this. We also appeal to base change for GL$(n)$ (\cite{dr:AC}). The proof given here is a strengthened form of the one in \cite{dr:Ra1}, and we avoid making use of either regularity or cuspidality (of $\pi$). This article is completely self-contained, however, and one does not need to refer to \cite{dr:Ra1}.

When the semisimplification of an $p$-adic representation $\tau$ decomposes as $\oplus_{j=1}^r \, tau_{j}$ with each $\tau_{j}$ irreducible of dimension $n_j$, we will say that~$\tau$ has decomposition type $(n_1, n_2, \dots, n_r)$. Similarly, if $\pi$ is an isobaric automorphic representation of GL$_n({\mathbb A}_F)$ of the form $\boxplus_{j=1}^r \, \pi_j$, with each $\pi_j$ cuspidal of GL$_{n_j}({\mathbb A}_F)$, we will say that $\pi$ has isobaric type $(n_1, n_2, \dots, n_r)$.

One motivation for the assertion of Theorem A (in \cite{dr:Ra1}) was to help establish the {\it irreducibility of $\rho$ for $\pi$ cuspidal and regular,algebraic} on GL$(4)/F$, at least {\it for large enough $p$}, generalizing the classical results of Ribet for $n=2$, and  Blasius-Rogawski (\cite{dr:BRo2}) for $n=3$ {\it and}
$\pi$ essentially self-dual:

\medskip

\noindent{\bf Theorem B} \, \it Let $F$ be totally real, $n \leq 4$, and $\pi$ an isobaric, algebraic automorphic representation of GL$_n({\mathbb A}_F)$ which is regular at every archimedean place, with associated Hodge-Tate representation $\rho$ of ${\mathfrak G}_F$. Assume that $\rho$ is crystalline and $p-1$ is greater than the twice the largest difference in the Hodge-Tate weights.
Then the (isobaric) type of $\pi$ determines the (reducibility) type of $\rho^{\rm ss}$. In particular, if $\pi$ is cuspidal, then $\rho$ is irreducible.
\rm

\medskip

It should be remarked that recently, such an irreducibility result for essentially selfdual, cuspidal $\pi$ has been established by F.~Calegari and T.~Gee (\cite{dr:CaG}) {\it with no condition on $p$} (and even for $n=5$), and the route they take appeals to, and further extends, recent modularity results, as well as \cite{dr:BCh} on signs of selfdual representations and the results on even Galois representations in \cite{dr:Ca}; see also \cite{dr:DV} and all other related references in \cite{dr:CaG}. In the non-selfdual case, the authors of \cite{dr:CaG} give a condensed form of a simple argument of \cite{dr:Ra1}, which we give in its entirety in section~3 below, both to correct an identity slightly and also to check that the isobaric type of $\pi$ agrees with the decomposition type of $\rho$.

Our proof of Theorem B makes use of Theorem A in conjunction with a theorem of Richard Taylor establishing the potential automorphy (over a totally real extension) of a class of $2$-dimensional $p$-adic representations, as well as some tricks involving $L$-functions. Our original approach was to require potential automorphy over a common totally real extension {\it simultaneously} for two $2$-dimensional odd representations of the type considered in \cite{dr:Ta2,dr:Ta3}. The proof given here is simpler and does not require this more stringent condition.  We also make use of the recent work \cite{dr:HLTT} on non-selfdual representations. In general, the arguments here are a bit more involved than in \cite{dr:Ra1} so as to deal (for Theorem A) also with the quasi-regular, non-regular case. If $F={\mathbb Q}$, in many cases one can also use \cite{dr:Ki} for the proof of the $(2,2)$-case of Theorem B (instead of \cite{dr:Ta2,dr:Ta3}), once the even summands are ruled out (as in Theorem A).

As a companion to Theorem B, we will also establish the following

\medskip

\noindent{\bf Theorem C} \, \it Let $F$ be a totally real number field, $p$ a rational prime. and $\sigma, \sigma'$ $2$-dimensional semisimple, odd, crystalline  $p$-adic representations of ${\mathfrak G}_F$ of the same weight $w$, such that each has distinct Hodge-Tate weights in an interval of length at most $(p-1)/2$. Assume moreover that $\sigma\oplus\sigma'$ is associated to an isobaric automorphic representation $\pi$ of GL$_4({\mathbb A}_F)$. Then there are a finite set $S$ of places of ${\mathbb Q}$ such that $L^S(1+w/2, \sigma) \, \ne \, 0$. Moreover,
$$
-{\rm ord}_{s=1} \, L^S(s, \sigma^\vee\otimes\sigma') \, = \, {\rm dim}(\sigma^\vee\otimes\sigma')^{{\mathfrak G}_F}.
$$
\rm

\medskip

This assertion is as predicted by the Tate conjecture when $\sigma, \sigma'$ occur in the $p$-adic \'etale cohomology of smooth projective varieties. It is essentially a consequence of Theorem B if we also assume that $\pi$ is regular and algebraic (in which case the oddness of $\sigma, \sigma'$ would in fact be a consequence, thanks to Theorem A). Furthermore, this result is stated (for $F={\mathbb Q}$) in \cite{dr:Ra1}, and referred to by Skinner and Urban in \cite{dr:SU} (Thm. 3.2.4), without even assuming that $\sigma \oplus \sigma'$ is automorphic, and this stronger version will be the subject of a sequel to this paper. However, Theorem C as stated above appears to be sufficient for the application in \cite{dr:SU}, section 3.2.6. To elaborate, in \cite{dr:SU}, $\sigma\oplus \sigma'$ is modular, corresponding to a cusp form $\Pi$ on GSp$(4)/{\mathbb Q}$. By Arthur's recent work (\cite{dr:A}), one can transfer $\Pi$ to an isobaric automorphic form $\pi$ on GL$(4)/{\mathbb Q}$ such that the degree $4$ $L$-functions of $\pi$ and of $\Pi$ coincide at all but a finite number of primes. (One can say more but it is not needed for this application.) Now Theorem C above applies with this $\pi$ associated to $(\sigma, \sigma')$. It will be interesting to extend this result without making the regularity assumption on the Hodge-Tate weights, possibly by making use of recent works of Calegari and Geraghty.

We would like to thank Don Blasius, Dipendra Prasad, Freydoon Shahidi, Chris Skinner and Eric Urban
for their interest and encouragement, Clozel for some remarks, and Richard Taylor
for helpful comments on an earlier ($\sim$ 2004) version of the preprint \cite{dr:Ra1}, which had clearly been inspired by \cite{dr:Ta2,dr:Ta3}. Thanks are also due to the referee for carefully reading the article and making helpful comments which improved the presentation. Finally, we happily acknowledge partial support from the NSF through the grant DMS-1001916.

\medskip

\section{Reductions}

\medskip

From here on, let $\pi$ be as in Theorem A, with associated $\rho$. Before deriving some lemmas, let us make some simple observations. To begin, note that
if~$\pi$ is an isobaric sum of idele class characters, then $\rho^{\rm ss}$ must be a direct sum of  one-dimensional $p$-adic representations. Indeed, each idele class character~$\nu$ appearing as an isobaric summand of $\pi$ is algebraic and corresponds to an $p$-adic character $\nu$ by \cite{dr:Se}, and the fact that $\rho^{\rm ss}$ is a direct sum of such characters follows by Tchebotarev. So we may take $2\leq n \leq 4$; of course the case $n=1$ is classical.

Next suppose there is a cuspidal automorphic representation $\beta$ of\linebreak GL$_2({\mathbb A}_F)$ occurring as an isobaric summand of $\pi$. Then by our hypothesis (in Theorem A), $\beta$ is a totally odd, algebraic cuspidal automorphic representation of GL$_2({\mathbb A}_F)$, which corresponds to a holomorphic Hilbert modular newform with multi-weight $(k_1, k_2, \dots, k_m)$, $m=[F:{\mathbb Q}]$, such that the $k_j\geq 1$ are all of the same parity. By Blasius-Rogawski (\cite{dr:BRo1}), Taylor (\cite{dr:Ta1}) in the regular case ($k_j \geq 2)$, and by Jarvis (\cite{dr:Jar}) when one of the $k_j$ is $1$, we know the existence of an irreducible $2$-dimensional $p$-adic representation $\tau$ of the absolute Galois group of $F$ attached to $\beta$, i.e., with coincident local factors almost everywhere. Then $\tau$ is totally odd, and there is nothing else to say if $n=2$; therefore we assume in the rest of this section that $n$ is $3$ or $4$. If we write $\pi=\beta \boxplus \beta'$, with $\beta'$ is an algebraic, isobaric
automorphic representation of GL$_{n-2}({\mathbb A}_F)$, then $\beta'$ also corresponds to an $(n-2)$-dimensional $p$-adic representation $\tau'$ which is totally odd if $n=4$.

Suppose $n=3$ and $\pi$ is a cuspidal. Then $\pi':=\pi \boxplus 1$ is an isobaric automorphic representation of GL$_4({\mathbb A}_F)$ satisfying the same hypotheses of Theorem A, with corresponding Galois representation $\rho' = \rho \oplus 1$. The proof of Theorem A for $\pi'$ will imply the same for $\pi$.

In view of the above remarks, we may assume from here on that $n=4$ and that either $\pi$ is cuspidal or isobaric of type $(3,1)$. Let $\Lambda^2(\pi)$ denote the isobaric automorphic representation of GL$_6({\mathbb A}_F)$ defined by Kim (\cite{dr:K}), which corresponds to $\Lambda^2(\rho)$.

Assume moreover that
$$
\rho^{\rm ss} \, \simeq \, \tau_{1} \, \oplus \, \tau_{2}, \, \, \, {\rm with} \, \, \, {\rm dim}(\tau_{j})=2, \, \, j=1,2.\label{dr:eqn1.1}
$$

\medskip

\noindent{\bf Lemma 1.2}\label{dr:lem1.2} \, \it The decomposition \eqref{dr:eqn1.1} of $\rho^{\rm ss}$
precludes the possibility of  $\pi$ being of type $(3,1)$, i.e.,
we cannot have an isobaric sum decomposition of the form
$$
\pi \, \simeq \, \eta \boxplus \nu,\leqno(\ast)
$$
with $\eta$ a cuspidal automorphic representation of GL$_3({\mathbb A}_F)$ and $\nu$ an idele class character. In fact, when $\pi$ is of type $(3,1)$, $\rho^{\rm ss}$ must be of the same type, i.e., $\rho$ must contain a sub or a quotient which is irreducible of dimension $3$.
\rm

\medskip

\begin{proof} Let $\omega$, resp. $\omega_\eta$, be the central character of $\pi$, resp. $\nu$, so that
$\omega_\nu=\omega\nu^{-1}$. Note that $(\ast)$ implies
$$
\Lambda^2(\pi) \, \simeq \, (\eta^\vee\otimes\omega_\pi\nu^{-1}) \boxplus (\eta\otimes \nu),
$$
which can be seen by checking the unramified local factors and applying the strong multiplicity one theorem (\cite{dr:JS1}). It follows that
\begin{equation}
\Lambda^2(\pi)\boxplus \nu^2 \boxplus \omega\nu^{-1} \, \simeq \, (\pi^\vee\otimes\omega\nu^{-1}) \boxplus (\pi\otimes \nu).\label{dr:eqn1.3}
\end{equation}
The algebraicity of $\pi$ implies the same for $\eta$, $\nu$ and $\omega$, and the latter pair defines (by \cite{dr:Se}) $p$-adic characters, again denoted by $\nu, \omega$, of the absolute Galois group of $F$. We get (by Tchebotarev)
\begin{equation}
\Lambda^2(\rho^{\rm ss}) \oplus \nu^2 \oplus \omega\nu^{-1} \, \simeq \, ({\rho}^{\rm ss})^\vee\otimes\omega\nu^{-1} \, \boxplus \, (\rho^{\rm ss}\otimes \nu).\label{dr:eqn1.4}
\end{equation}
On the other hand, applying the exterior square operation to both sides of the decomposition \eqref{dr:eqn1.1} yields
\begin{equation}
\Lambda^2(\rho^{\rm ss}) \oplus \nu^2 \oplus \omega\nu^{-1} \, \simeq \, \left(\tau_{1}\otimes\tau_{2}\right) \oplus \omega_{1} \oplus \omega_{2} \oplus \nu^2 \oplus \omega\nu^{-1},\label{dr:eqn1.5}
\end{equation}
where $\omega_{j}$ is the determinant of $\tau_{j}$.
Comparing \eqref{dr:eqn1.4} and \eqref{dr:eqn1.5}, and using \eqref{dr:eqn1.1} and the isomorphism $\tau_{j}^\vee \simeq \tau_{j}\otimes \omega_{j}^{-1}$, we see that $\tau_{1}$ and $\tau_{2}$ must both be reducible. If we write
$$
\tau_{j} \, \simeq \, \mu_{j} \oplus \mu_{j}', \, \, {\rm dim}(\mu_{j}) = {\rm dim}(\mu'_{j}) = 1,
$$
then the characters on the right are both Hodge-Tate, since $\rho$ has that property (by assumption), hence are locally algebraic and correspond to idele class characters of $F$ (cf. \cite{dr:Se}). This then forces the isobaric decomposition of the form
$$
\pi \, \simeq \, \mu_1 \boxplus \mu_1' \boxplus \mu_2 \boxplus \mu_2',
$$
which contradicts the cuspidality of $\eta$.

To finish the proof of the Lemma, we still need to check that $\rho$ cannot be irreducible when $(\ast)$ holds. For this, note that \eqref{dr:eqn1.4} still holds because its proof did not use \eqref{dr:eqn1.1}. But then, since the left hand side contains $1$-dimensionals, $\rho$ cannot be irreducible.
\end{proof}

\medskip

\noindent{\bf Lemma 1.3}\label{dr:lem1.3} \it Let $\pi$ be cuspidal and satisfy the hypotheses of Theorem A, together
with the decomposition \eqref{dr:eqn1.1} of the semisimplification of the associated $\rho$.
Let $K$ be a totally real number field which is solvable and Galois over $F$. If the base change $\pi_K$ of $\pi$ to GL$(4)/K$ is Eisensteinian, then $\rho$  must be irreducible.
\rm

\medskip

\begin{proof} If $\pi_K$ is Eisensteinian, the cuspidality of
$\pi$ implies that Gal$(K/F)$ acts transitively on the set of
cuspidal automorphic representations occurring in the isobaric sum
decomposition of $\pi_K$ (cf. \cite{dr:AC}). This forces a decomposition $\pi_K
\simeq \beta_1 \boxplus \beta_2$, where each $\beta_j$ is an
isobaric automorphic representation of GL$(2, {\mathbb A}_K)$. It follows
that there are intermediate fields $L, E$ such that $F \subset E
\subset L \subset K$ with $[L:E]=2$, and a cuspidal automorphic
representation $\beta$ of GL$(2, {\mathbb A}_{L})$ such that $\pi_E$ is
cuspidal and is the automorphic induction $I_{L}^E(\beta)$.
Now since $\pi_E$ is the base change of $\pi$, at every archimedean
place $w$ of $E$, the eigenvalue set of complex conjugation in the
associated $4$-dimensional representation $\sigma_w$ of $W_{E_w}\simeq W_{\mathbb R}$
must remain $\{1,-1,1,-1\}$. As $L$ is totally real, $w$ splits in it,
say into $\{\tilde{w}_1, \tilde{w}_2\}$, and so
$\beta_{\tilde{w}_j}$ is forced to be odd for each $j$. In other words, $\beta$ is a totally odd,
algebraic cusp form on GL$(2)/L$, and by \cite{dr:BRo1}, \cite{dr:Ta1} and \cite{dr:Jar}, it corresponds to a
$2$-dimensional, irreducible, odd Galois representation $r$. Then the induction $R$, say,
of $r$ to the absolute Galois group ${\mathfrak G}_E$ of $E$
corresponds by functoriality to $\pi_E$. Moreover, the cuspidality of $\pi$ implies that $\beta$,
and hence $r$, is not $\theta$-invariant, where $\theta$ denotes the non-trivial automorphism of $L/E$.
Then $R$ is irreducible by Mackey, and must be the restriction of $\rho$ to ${\mathfrak G}_E$; this can be checked, for example, at almost all places and then deduced by Tchebotarev. In any case, it implies that $\rho$ itself must be irreducible.
\qed\end{proof}

Note that in this section we did not need $\pi$ to be quasi-regular, only that it is totally odd.

\section{Proof of Theorem A}

Fix a totally real number field $F$.
In view of the Lemmas of the previous section, we may assume that $n=4$ with $\pi$ satisfying the hypotheses
of Theorem A with associated $\rho$, and that $\pi$ remains
cuspidal upon base change to any finite solvable Galois extension $K$ which is totally real.

Suppose we have a decomposition of the form \eqref{dr:eqn1.1}, i.e., with $\rho^{\rm ss}$ being the direct sum of two $2$-dimensional $p$-adic representations $\tau_{1}$, $\tau_{2}$ of ${\mathfrak G}_F$, with $\omega_{j}={\rm det}(\tau_{j})$, which we will also view (by class field theory) as
idele class characters $\omega_j$ of $F$.
The Hodge-Tate hypothesis on $\rho$ implies that as a $p$-adic character,
$\omega_{j}$ is locally algebraic for $j =1,2$, and is thus
associated to an algebraic Hecke character, again denoted $\omega_j$, of the idele
class group $C_F$ of (the totally real) $F$, which must be an integral power of the norm
character times a finite order character $\nu_j$. It follows that as a $p$-adic character,
$$
\omega_{j} \, = \chi^{a_j}\nu_j,
$$
where $\chi$ is the $p$-adic cyclotomic character.
Clearly, $\omega={\rm det}(\rho) =
\omega_{1}\omega_{2}$. Note that $\tau_j$ is even relative to some $c$ iff $a_j$ and $\nu_j$ have
the same parity at an archimedean place.

\medskip

\noindent{\bf Lemma 2.1}\label{dr:lem2.1} We have $a_1 = a_2$. Also, $\nu_1$
and $\nu_2$ have the same parity.
\rm

\medskip

\begin{proof} We may, and we will, assume that $a_1 \geq a_2$. A
twist of $\pi$ by $\vert \cdot \vert^t$ is, for some $t \in {\mathbb R}$,
unitary; call this representation $\pi^u$. Since $\pi$ is
algebraic, $t$ lies in $\frac12{\mathbb Z}$. Then at any finite place $v$ where
$\pi$ is unramified, we know by the Rankin-Selberg theory that the
inverse roots $\alpha_{j,v}$, $1 \leq j \leq 4$, defining the
local factor of $\pi^u$, satisfy $\vert \alpha_{j,v}\vert <
(Nv)^{1/2}$.
Consequently, any inverse root of $\Lambda^2(\pi^u_v)$ is strictly
bounded in absolute value by $Nv$. Now,
$$
\Lambda^2(\pi^u) \, \simeq \, \Lambda^2(\pi \otimes \vert \cdot
\vert^t) \, \simeq \, \Lambda^2(\pi) \otimes \vert \cdot
\vert^{2t}.
$$
So $\vert a_j - 2t\vert < 1$ for $j=1,2$, and since $a_j, 2t \in
{\mathbb Z}$, $a_1=2t=a_2$ (as desired).

Consequently, the central character of $\pi$, which is associated to $\omega={\rm det}(\rho)$, equals
$\nu_1\nu_2\vert \cdot\vert^{2a}$, where $a=a_1=a_2$. Since $\omega$ is even, $\nu_1$ and $\nu_2$ have the same parity.
\qed\end{proof}

\noindent {\bf Assume} the following:
\begin{equation}\label{dr:eqn2.1}
\text{For $j=1,2$, $\nu_j$ has the same parity as $a$ at some archimedean
place $v_0$.}
\end{equation}

\bigskip

We will obtain a contradiction below, leading to Theorem A.

\medskip

\noindent{\bf Lemma 2.2}\label{dr:lem2.2} \, \it Let $K$ be a finite abelian
extension of $F$ in which $v_0$ splits. Then $\omega_{j,K}$ ($=\omega_j\circ N_{K/F}$)
cannot occur, for either $j$,
in the isobaric sum decomposition of $\Lambda^2(\pi_K)$. In
particular, $L(s, \Lambda^2(\pi)\otimes\omega_{j,K}^{-1})$ has no
pole at $s=1$.
\rm

\medskip

It should be noted that this crucial Lemma will be false if
$\omega_j$ were to be totally odd; in that case $L(s, \pi; {\rm
sym}^2\otimes\omega_{j,K}^{-1})$ will be regular at $s=1$.

\medskip

\begin{proof}  First consider when $\pi$ is regular at the archimedean place $v_0$,
so that the associated $4$-dimensional representation of $W_{\mathbb R}$ is of the form
$$
\sigma_{v_0} \, \simeq \, I(\xi) \oplus I(\xi'),
$$
for distinct, non-conjugate characters $\xi:z\mapsto z^{p}\overline{z}^{w-p}$ and
$\xi':z\mapsto z^{r}\overline{z}^{w-r}$ of ${\mathbb C}^\ast$, where $I$ denotes the
induction from ${\mathbb C}^\ast$ to $W_{\mathbb R}$. Then
$I(\xi) \otimes I(\xi')$ decomposes as $I(\xi\xi') \oplus
I(\xi{\overline \xi'})$. It follows that
\begin{equation}
\Lambda^2(\sigma_\infty) \, \simeq \, I(z^{p+r}{\overline
z}^{2w-p-r}) \oplus I(z^{p+w-r}{\overline z}^{w-p+r}) \oplus {\it
sgn}^{w+1}\vert\cdot\vert^w \oplus {\it
sgn}^{w+1}\vert\cdot\vert^w,\label{dr:eqn2.3a}
\end{equation}
where $w$ is the archimedean weight of $\pi$ (with $w-a_j\in {\mathbb Z}$).
Since $K/F$ is abelian, the base change
$\pi_K$ (of $\pi)$ makes sense (\cite{dr:AC}), and for any archimedean
place $\tilde{v_0}$ of
$K$ above $v_0$, $K_{\tilde{v_0}} \simeq {\mathbb R}$ (since $v_0$
splits in $K$ by assumption)
and $\sigma_{\tilde{v_0}} \simeq \sigma_{v_0}$ is the
parameter of $\pi_{K,\tilde{v_0}}$ if $\tilde{v_0} \mid {v_0}$.
So the occurrence of $\omega_{j,K}$ (for either
$j$) in the isobaric decomposition of $\Lambda^2(\pi_K)$ implies
that $\omega_{j,K,\tilde{v_0}}$ will occur in $\Lambda^2(\sigma_\infty)$, and
so must equal ${\it sgn}^{w+1}\vert\cdot\vert^w$. Now since
$\omega_{j,v_0}$ is of the form $\omega_{j,v_0}\vert\cdot\vert^a$, we see that
$w = a$ and that the parity of $\omega_{j}$ at $v_0$ is the opposite of that
of $a$, a contradiction!

Next look at when $\pi_{v_0}$ is quasi-regular, but not regular. Then we must have
$$
\sigma_{v_0} \, \simeq \, I(\xi) \oplus (1+{\operatorname{sgn}})\vert \cdot \vert^{w/2},
$$
where again $\xi$ is given as above. In this case, since $I(\xi)$ is twist invariant under
${\operatorname{sgn}}$, the sign character, we obtain
\begin{equation}
\Lambda^2(\sigma_\infty) \, \simeq \, \left(I(\xi) \oplus {\operatorname{sgn}}\vert \cdot \vert^{w}\right)^{\oplus 2}.
\label{dr:eqn2.3b}
\end{equation}
Again, we see that $\omega_{j,\tilde{v_0}}$ cannot be a summand on the right because of parity.
\qed\end{proof}

As the proof shows, $a$ is necessarily the archimedean weight of $\pi$.

\noindent{}Let
$$
\alpha = \nu_1/\nu_2, \, \, \nu=\nu_1\nu_2, \quad {\rm and} \quad E = F(\alpha, \nu),
$$
the compositum of the cyclic extensions of $F$ cut out by $\alpha$ and $\nu$.
Then by Lemma \ref{dr:lem2.1}, $\alpha=\omega_1/\omega_2$, and $\omega=\nu\vert\cdot\vert^{2a}$.
Moreover, the total evenness of
$\alpha$ and $\nu$ implies that $E$ is totally real. By construction,
$\nu_{1,E} = \nu_{2,E}$ and $\nu_E =1$, so if we put
$$
\mu: = \nu_{1,E}\vert\cdot\vert_E^a,
$$
then $\mu^2$ is $\vert\cdot\vert_E^{2a}$ on $C_E$.

\newpage

The reason for considering $E$ is that $\pi$ becomes essentially selfdual over~$E$:

\medskip

\noindent{\bf Lemma 2.4}\label{dr:lem2.4} \, \it We have
\begin{enumerate}
\setlength{\itemsep}{-0.4ex}
\item[{\rm (a)}] $\pi_E^\vee \, \simeq \, \pi_E \otimes \mu^{-1}$;
\item[{\rm (b)}] The incomplete $L$-function $L^T(s, \pi_E; {\rm sym}^2
\otimes \mu^{-1}))$ has a pole at $s=1$ of order $1$, where $T$ is
a finite set of places of $E$ containing the archimedean and
ramified primes;
\item[{\rm (c)}] $\Lambda^2(\pi_E\otimes\vert\cdot\vert_E^{-a})$ is selfdual.
\end{enumerate}
\rm

\medskip

\begin{proof}  (a) \, By assumption, $\rho^{\rm ss}$ is
$\tau_{1} \oplus \tau_{2}$, and for either $j$,
$\tau_{j}^\vee \simeq \tau_{j} \otimes
\omega_{j}^{-1}$. So with $\mu$ being the restriction of
$\omega_{j}$ to $E$, then
$$
{\rho_{E}^{\rm ss}}^\vee \, \simeq \, \rho_{E}^{\rm ss}
\otimes \mu^{-1}.
$$
Comparing $L$-functions, we get
$$
L^T(s, \pi_E^\vee) \, = \, L^T(s, \pi_E \otimes \mu^{-1}).
$$
It follows by the strong multiplicity one theorem for isobaric
automorphic representations (\cite{dr:JS1}) that
$\pi_E^\vee$ is isomorphic to $\pi_E \otimes \mu^{-1}$.

\medskip

(b) \, We have
$$
L^T(s, \pi_E \times \pi_E \otimes \mu^{-1}) \, = \, L^T(s, \pi_E,
\Lambda^2 \otimes \mu^{-1}))L^T(s, \pi_E, {\rm sym}^2 \otimes
\mu^{-1})).
$$
One knows by the Rankin-Selberg theory that the $L$-function on
the left has a pole at $s=1$ of order $\geq 1$, and so the assertion
follows in view of Lemma \ref{dr:lem2.2} (with $K=E$).

\medskip

(c) \, This follows from (a).
\qed\end{proof}

For any isobaric automorphic representation $\Pi$ of GL$(6)/F$ and
a character $\xi$ of $F$, let $L^T(s, \Pi; \Lambda^3 \otimes \xi)$
denote, for a finite set $T$ of places, the incomplete
$\xi$-twisted {\it exterior cube} $L$-function of $\Pi$ of degree
$20$. One knows --- see \cite{dr:Sh1}, Corollary 6.8, that this
$L$-function admits a meromorphic continuation to the whole
$s$-plane and satisfies a standard functional equation.

\medskip

\noindent{\bf Proposition 2.5} \label{dr:prop2.5} \, \it Let $\Pi=
\Lambda^2(\pi)$, and $T$ a sufficiently large finite set of
places of $E$ containing the archimedean and ramified places as well as
those dividing $p$. Then for any character $\xi$ of $E$,
\begin{enumerate}
\setlength{\itemsep}{-0.4ex}
\item[{\rm (a)}]$L^T(s, \Pi_E; \Lambda^3\otimes \xi) = L^T(s, \pi_E; {\rm
sym}^2 \otimes \xi)^2$;
\item[{\rm (b)}]$L^T(s, \Pi_E;
\Lambda^2\otimes \xi)L^T(s, \xi) = L^T(s, \pi_E \times \pi_E
\otimes \xi\mu)$.
\end{enumerate}
\rm

\newpage

\noindent{\bf Corollary 2.6}\label{dr:cor2.6} \rm \hfill
\begin{enumerate}
\setlength{\itemsep}{-0.4ex}
\item[{\rm (a)}]$-{\rm ord}_{s=1} L^T(s,\Pi_E; \Lambda^3\otimes \mu) \geq 2$;
\item[{\rm (b)}]${\rm ord}_{s=1}L^T(s,\Pi_E; \Lambda^2\otimes\mu^{-2}) =
0$.
\end{enumerate}
\rm

\medskip

Clearly, part (a) (resp. (b)) of the Corollary follows from part
(b) (resp. (a)) of Lemma \ref{dr:lem2.4} and part (a) (resp. (b) of
Proposition \ref{dr:prop2.5}.

\medskip

Clearly, part (a) (resp. (b)) of the Corollary follows from part
(b) (resp. (a)) of Lemma 2.4 and part (a) (resp. (b)) of
Proposition 2.5.

\qed

\medskip

To prove the Proposition we need the following

\medskip

\noindent{\bf Lemma 2.7}\label{dr:lem2.7} \, \it Let $\sigma = \sigma_1 \oplus
\sigma_2$ be a representation of a group $G$, with each $\sigma_j$
being of dimension $2$ with determinant
$\alpha_j$. Then ${\rm sym}^2(\sigma)$ is $\sigma_1 \otimes \sigma_2 \oplus
{\rm sym}^2(\sigma_1) \oplus {\rm sym}^2(\sigma_2)$. Moreover,
\begin{multline*}
\Lambda^3(\Lambda^2(\sigma)) \, \simeq \\
\left(\sigma_1\otimes\sigma_2\otimes\alpha_1\alpha_2\right)^{\oplus
2} \, \oplus \, {\rm
sym}^2(\sigma_1)\otimes(\alpha_1\alpha_2\oplus\alpha_2^2) \,
\oplus \, {\rm
sym}^2(\sigma_2)\otimes(\alpha_1\alpha_2\oplus\alpha_1^2),\\
\sigma^{\otimes 2} \, \simeq {\rm sym}^2(\sigma_1)\otimes
\alpha_1 \oplus {\rm sym}^2(\sigma_2)\otimes \alpha_2 \oplus
(\sigma_1\otimes\sigma_2)^{\oplus 2},
\end{multline*}
and
$$
\Lambda^2(\Lambda^2(\sigma)) \, \simeq \,
(\sigma_1\otimes \sigma_2)^{\oplus 2} \, \oplus {\rm sym}^2(\sigma_1)\otimes\alpha_2 \,
\oplus \, {\rm sym}^2(\sigma_2)\oplus \alpha_1.
$$
In particular, if $\alpha_1=\alpha_2=\alpha$, we have
\begin{enumerate}
\item[{\rm (a)}] $\Lambda^3(\Lambda^2(\sigma))\otimes\alpha^{-2} \, \simeq \, {\rm
sym}^2(\sigma)^{\oplus 2}$ 

\medskip

and

\medskip

\item[{\rm (b)}]
$\Lambda^2(\Lambda^2(\sigma))\oplus \alpha^2 \, \simeq \, \sigma
\otimes \sigma\otimes\alpha$.
\end{enumerate}
\rm

\medskip

\noindent{\it Proof} \, The first identity (involving the symmetric squa\-re)
is evident. For the
second, observe that since $\Lambda^2(\sigma)$ is $\sigma_1
\otimes\sigma_2\oplus\alpha_1\oplus\alpha_2$,
\begin{enumerate}
\item[(i)] $\Lambda^3(\Lambda^2(\sigma)) \, \simeq
\Lambda^3(\sigma_1\otimes\sigma_2) \oplus
\Lambda^2(\sigma_1\otimes\sigma_2)\otimes\left(\alpha_1\oplus\alpha_2\right)
\oplus\sigma_1\otimes\sigma_2\otimes\alpha_1\alpha_2$.

\medskip

We have

\medskip

\item[(ii)] $\Lambda^2(\sigma_1\otimes\sigma_2) \, \simeq \, {\rm
sym}^2(\sigma_1)\otimes\Lambda^2(\sigma_2)
\oplus\Lambda^2(\sigma_1)\otimes{\rm sym}^2(\sigma_2)$.
\end{enumerate}
Moreover, the non-degenerate $G$-pairing
$$
\sigma_1\otimes\sigma_2 \times \Lambda^3(\sigma_1\otimes\sigma_2)
\, \to \, \Lambda^4(\sigma_1\otimes\sigma_2) \, = \, {\rm
det}(\sigma_1\otimes\sigma_2) = \alpha_1^2\alpha_2^2
$$
identifies $\Lambda^3(\sigma_1\otimes\sigma_2)$ with
$(\sigma_1\otimes\sigma_2)^\vee \otimes \alpha_1^2\alpha_2^2$,
which is isomorphic to $\sigma_1\otimes\sigma_2 \otimes
\alpha_1\alpha_2$. The second identity and part (a) of the Lemma
now follow by putting these together.

Moreover,
$$
\Lambda^2(\Lambda^2(\sigma)) \, \simeq \,
\Lambda^2(\sigma_1\otimes\sigma_2) \oplus \alpha_1\alpha_2 \oplus
\left(\sigma_1\otimes\sigma_2\otimes(\alpha_1\oplus\alpha_2)\right).
$$
Applying (ii), we also get the third identity and part (b).

\qed

\noindent{\bf Proof of Proposition 2.5} \, Since $\rho$ is associated
to $\pi$, there is a finite set $S$ of places of
$F$ containing the ramified places and $\infty$, such that for
all $u \notin S$,
$$
L(s, \pi_u) \, = \, L(s, \rho_{u}^{\rm ss}),
$$
which is just an equality of $4$-tuples of inverse roots defining
the respective $L$-factors. It follows using \eqref{dr:eqn1.4} that for any
such $u$,
$$
L(s, \Lambda^2(\pi)_u) \, = \, L(s, \tau_{1,u} \otimes
\tau_{2,u})L(s,\omega_{1,u})L(s,\omega_{2,u}).
$$
Since by construction, $\omega_{1}$ and $\omega_{2}$ become
the same over $E$, it follows that at any place $v$ of $E$
outside the inverse image $T$ of $S$, we can write
$$
\sigma_v(\pi_E) \, \simeq \, \sigma_{1,v} \oplus \sigma_{2,v},
$$
with $\sigma_{1,v}$, $\sigma_{2,v}$ both being $2$-dimensional of
the same determinant. Here, $\sigma_v(\pi_E)$ is the
$4$-dimensional representation of the Weil group of $F_{1,v}$ attached
to $\pi_{E,v}$. Now we are done by appealing to Lemma \ref{dr:lem2.7}.

\qed

\medskip

\noindent{\bf Proof of Theorem A (contd.)} \hfill

As noted at the beginning of this section, we have already reduced to the case when $\pi$ remains cuspidal over any solvable totally real extension of $F$, in particular over $E$. Recall that $\pi_E$ is essentially selfdual relative to the character $\mu$, and we may replace $\pi$ by its twist by a power of $\vert\cdot\vert$ so that $\Pi_E$ is selfdual.

Suppose $\Pi_E=
\Lambda^2(\pi_E)$ is cuspidal. Then by a result of Ginzburg and
Rallis (cf. \cite{dr:GR}, Theorem 3.2) we know that $L^T(s, \Pi_E;
\Lambda^3\otimes \xi)$, for any character $\xi$,
will have at most a simple pole at $s=1$,
for any finite set $T$ of places of $E$ containing the archimedean
and ramified places. (If $\Pi_E$ has a supercuspidal component, this
also follows from \cite{dr:KSh}, where one finds such a result even for
the complete $L$-function.) But this contradicts Corollary~\ref{dr:cor2.6}, part (a).

So $\Pi_E$ is not cuspidal, and we may write
\begin{equation}
\Pi_E \, \simeq \, \boxplus_{j=1}^m \beta_j, \, m > 1,\label{dr:eqn2.11}
\end{equation}
where each $\beta_j$ is a cuspidal automorphic representation of
GL$(n_j, {\mathbb A}_E)$ with $\sum_j n_j = 6$ and $n_i \leq n_j$ if $i
\leq j$. Note that each $\beta_j$ will necessarily be algebraic.

It remains to get a contradiction when $\Pi_E$ is Eisenteinian. Note that since
$\pi_E$ has trivial central character, $\Lambda^2(\pi_E)$ is selfdual.

\newpage

\noindent{\bf Lemma 2.8}\label{dr:lem2.8} \, \it
$n_1 \geq 2$ and $m\leq 3$.
\rm

\medskip

\noindent{\it Proof}. \, First note that if some $n_j=1$, then twisting by the
inverse of the corresponding character $\beta_j$,
we get a pole at $s=1$ of the $L$-function
$$
L^T(s, \pi_E \times \pi_E\otimes \beta_j^{-1}) \, = \, L^T(s, \Pi_E\otimes\beta_j^{-1})L^T(s,\pi_E, {\rm sym}^2\otimes
\beta_j^{-1}),
$$
because $\prod_{j\geq 2} L^T(s,\beta_j)$ and $L(s, \pi_E; {\rm sym}^2)$ has no zero at $s=1$ (cf. \cite{dr:Sh1}). Consequently, by the cuspidality of $\pi_E$ and \cite{dr:JS1}, we obtain
$$
\pi_E^\vee \, \simeq \, \pi_E\otimes\beta_j^{-1} \, \, \, ({\rm if} \, \, \, n_j=1).
$$
Moreover, at any archimedean place $v$ of $E$ lying above above $u$ (of $F$), the expressions \eqref{dr:eqn2.3a} and
\eqref{dr:eqn2.3b} (for the possible shapes of) $\Lambda^2(\sigma_{u}(\pi))$
imply, since $u$ splits in $E$, that {\it at most two of the
$n_j$ could be $1$}, and the corresponding $\beta_j$ must be totally odd. Suppose $n_1, n_2 =1$. Since $\pi$ is cuspidal,
$$
\gamma: = \beta_1\beta_2^{-1} \, \ne 1,
$$
and $\pi$ admits a self-twist by this totally even character,
$\pi$ must be
induced from the totally real quadratic extension $L=E(\gamma)$ of $E$; so $\pi_L$ is not cuspidal. But this cannot happen, as we have earlier reduced (see the first paragraph of this section) to the situation where the base change of $\pi$ to any finite solvable Galois, totally real extension remains cuspidal. Hence $n_2 \geq 2$ and $m \leq 3$.

If $n_1=1$, then there are then two possibilities for the type of $\Pi_E$, namely $(1,5)$ and $(1,2,3)$. In the former case,
$$
L^T(s, \Pi_E; \Lambda^3) \, = \, L^T(s, \beta_2; \Lambda^2 \otimes
\beta_1)L^T(s, \beta_2; \Lambda^3).
$$
The second $L$-function can be identified with $L^T(s,
\Lambda^2(\beta_2)^\vee \otimes \lambda)$, where $\lambda$ is the central
character of $\beta_2$. So each factor on the right of (5.13) is
an abelian twist of the exterior square $L$-function of a cusp
form on GL$(5)/F$, and by a theorem of Jacquet and Shalika
(\cite{dr:JS2}), it admits no pole at $s=1$ (because $5$ is odd). This
contradicts part (a) of Corollary \ref{dr:cor2.6}, and so the case $(1,5)$ cannot happen.

Suppose $(n_1,n_2,n_3)=(1,2,3)$. Here each $\beta_j$ must be selfdual. In particular, the square of the idele class character $\beta_1$ is $1$; similarly for the square of the central characters $\delta_j$ of $\beta_j$ for $j=2,3$. But $\beta_1$ must be non-trivial, because the $L$-function of $\Pi_E$ has no pole at $s=1$. Similarly, $L^T(s, \delta_2)$ divides $L^T(s, \Pi_E; \Lambda^2)$, which (by part (b) of Corollary \ref{dr:cor2.6} is invertible at $s=1$, implying that $\delta_2$ is non-trivial. Consequently, $\beta_2 \simeq \beta_2^\vee =\beta_2\otimes \delta_2$ and so $\beta_2$ is dihedral. Moreover, the selfduality of the GL$(3)$-cusp form $\beta_3$ implies (cf. \cite{dr:Ra5}, for example) that it is of the form sym$^2(\eta)\otimes \xi$ for a cusp form $\eta$ on GL$(2)/E$ and a character $\xi$.
A direct computation yields (using $\Lambda^2(\beta_3)\simeq \beta_3$ by selfduality):
\begin{multline*}
L^T(s, \Pi_E; \Lambda^3) =\\
L^T(s, \beta_1\times\delta_2)L^T(s,\beta_1\otimes\beta_2\times{\rm sym}^2(\eta)\otimes\xi)L^T(s, (\beta_1\delta_3^{-1})\otimes\beta_3)\times\\
L^T(s, \delta_2\otimes\beta_3)L(s,\xi_3).
\end{multline*}

Since the left hand side equals $L^T(s, \pi_E; {\rm sym}^2)^2$, the incomplete $L$-function
$L^T(s, \beta_1\xi\otimes\beta_2\times{\rm sym}^2(\eta))$ must be the square of a degree $3$ automorphic $L$-function $L(s)$. This forces $\eta$ to be dihedral, which is not possible since ${\rm sym}^2(\eta)$ is not cuspidal. Hence the $(1,2,3)$ case is not possible either.

\qed

\medskip

\subsection*{End of proof of Theorem A}

\medskip

To recap, we have already established Theorem A when $\Pi_E$ is cuspidal, and we may assume by Lemma \ref{dr:lem2.8} that $\Pi_E$ is an isobaric sum as in \eqref{dr:eqn2.11} with $n_1\geq 2$ and $m\leq 3$. As above, we may also replace $\pi$ by its twist by a power of $\vert\cdot\vert$ to assume that $\Pi_E$ is selfdual.

Suppose $\Pi_E$ is of {\it type $(2,2,2)$}, i.e., an isobaric sum of three algebraic cusp
forms $\beta_1, \beta_2, \beta_3$ on GL$(2)/E$, with $\alpha_j$ being the central character of $\beta_j$, we get (for suitable finite set $T$ of places containing the archimedean ones):
\begin{equation}
L^T(s, \Pi_E; \Lambda^3) \, = \, L^T(s, (\beta_1\boxtimes\beta_2)\times\beta_3)
\prod_{i\ne j} L^T(s, \beta_i\otimes\alpha_j).\label{dr:eqn2.13}
\end{equation}
We know by part (b) of Corollary \ref{dr:cor2.6} that this exterior cube $L$-function of $\Pi_E$ has a pole of order $\geq 2$ at $s=1$. Since each $L^T(s, \beta_i\otimes\alpha_j)$ is invertible at $s=1$ (as $\beta_i$ is cuspidal), we deduce from \eqref{dr:eqn2.13} that $L^T(s, (\beta_1\boxtimes\beta_2)\times\beta_3)$ must have at least a double pole, which implies (\cite{dr:Ra3}) that $\beta_1\boxtimes\beta_2$ must be of type $(2,2)$ and contain $\beta_2^\vee$ as an isobaric summand with multiplicity $2$. This forces (see \cite{dr:PRa}) each $\beta_i$ to be dihedral, in fact corresponding to a character $\chi_i$ of a common quadratic extension $E/E$. Then each $\beta_i$ is associated to an irreducible $2$-dimensional representation $\tau_i$ of ${\mathfrak G}_E$. It follows that the restriction of $\Lambda^2(\rho^{\rm ss})$ to ${\mathfrak G}_E$ must be $\tau_1 \oplus \tau_2 \oplus \tau_3$, which contradicts the fact that the (supposed) decomposition \eqref{dr:eqn1.1} results in $\Lambda^2(\rho^{\rm ss})$ having one-dimensional summands. Hence $\Pi_E$ cannot be of type $(2,2,2)$.

Next suppose $\Pi_E$ is of {\it type $(2,4)$}, i.e., with $\beta_1$, resp. $\beta_2$, being a selfdual cusp form on GL$(2)/E$, resp. GL$(4)/E$. Let $\xi_j$ be the central character of $\beta_j$ whose square is $1$. Then $L^T(s, \xi_1)$ divides $L^T(s, \Pi_E, \Lambda^2)$, which has no pole at $s=1$; thus $\xi_1\ne 1$, implying that the selfdual $\beta_1$ must be dihedral, say of the form $I_L^E(\chi)$ for a character $\chi$ of a quadratic extension $L$ of $E$. We obtain
\begin{equation}
L^T(s, \Pi_E; \Lambda^3) \, = \, L^T(s, \beta_1\times\Lambda^2(\beta_2))
L^T(s, \xi_1\otimes\beta_2)^2.\label{dr:eqn2.14} 
\end{equation}
The occurrence of at least a double pole then forces the isobaric decomposition
\begin{equation}
\Lambda^2(\beta_2) \, \simeq \, (\beta_1\otimes \xi_1)^{\boxplus 2}\boxplus I_L^E(\chi'),\label{dr:eqn2.15}
\end{equation}
with $I_L^E(\chi\chi')\simeq I_L^E(\chi^\theta\chi')$, where $\theta$ is the non-trivial automorphism of $L/E$. This forces $\chi'$ to be $\theta$-invariant, forcing $I_L^E(\chi')$ to be Eisensteinian and of the form $\lambda\boxplus\lambda\delta$, where $\lambda$ is the restriction of $\chi'$ to $C_E$ and $\delta$ is the quadratic character attached to $L$. Since $\beta_2$ is selfdual, it follows that it admits a self-twist by $\delta$ and is hence induced from $L$; say $\beta_2= I_L^E(\eta)$, for a cusp form $\eta$ on GL$(2)/L$ with central character $\gamma$. Then we have
\begin{equation}
\Lambda^2(\beta_2) \, \simeq \, {\rm As}_{L/E}(\eta)\otimes\delta \, \boxplus \, I_L^E(\gamma),\label{dr:eqn2.16}
\end{equation}
where As$_{L/E}(\eta)$ is the Asai transfer of $\eta$ to a form on GL$(4)/E$ (\cite{dr:Ra4}). Comparing \eqref{dr:eqn2.15} and \eqref{dr:eqn2.16} we see that this Asai representation is not cuspidal and of type $(2,2)$. Then by Theorem B of \cite{dr:PRa}, we see that $\eta$ must be dihedral, induced by a character $\lambda$, say, of a biquadratic extension $N$ of $E$ containing $L$. Thus
$$
\beta_2 \, \simeq \, I_L^E(\eta) \, \simeq \, I_N^E(\lambda),
$$
with $\lambda$ algebraic. So we may attach to $\beta_2$ an irreducible $4$-dimensional representation $\sigma$ obtained by Galois induction from the $p$-adic Galois character $\lambda$ associated to $\lambda$ in \cite{dr:Se}. Since $\Lambda^2(\pi_E)=\beta_1\boxplus \beta_2$, we obtain (by Tchebotarev) an isomorphism of ${\mathfrak G}_E$-modules
\begin{equation}
\Lambda^2(\rho^{\rm ss}) \, \simeq \, {\rm Ind}_L^E(\chi) \oplus {\rm Ind}_N^E(\lambda),\label{dr:eqn2.17}
\end{equation}
with the two summands on the right being irreducible. This contradicts the  fact that \eqref{dr:eqn1.1} implies that the left hand side of \eqref{dr:eqn2.17} admits one-dimensional summands, already over $F$. Hence the type $(2,4)$ is not possible in our case.

Finally, suppose $\Pi_E$ is of {\it type $(3,3)$}, and write
$$
\Lambda^2(\pi_E) \, \simeq \, \beta \boxplus \beta',
$$
with $\beta, \beta'$ cusp forms on GL$(3)/E$. Since $E$ is a totally real extension of $F$, the archimedean type is the same as over $F$. It follows from \eqref{dr:eqn2.3a}, \eqref{dr:eqn2.3b} that $\beta$, $\beta'$ must both be {\it regular} algebraic. Either they are both selfdual or duals of each other. In the former situation, one has known since \cite{dr:Pic} how to attach irreducible $3$-dimensional $p$-adic representations $\sigma$, $\sigma'$ to them; see also \cite{dr:Cl2}, \cite{dr:PL}. In the non-selfdual case, one still has such a construction in the recent work \cite{dr:HLTT}. Combining this with the exterior square of \eqref{dr:eqn1.1}, we get (by Tchebotarev) the following as ${\mathfrak G}_E$-modules:
$$
\tau\otimes\tau' \oplus \omega\oplus\omega' \, \simeq \, \Lambda^2(\rho^{\rm ss}) \, \simeq \, \sigma \oplus \sigma',
$$
which gives the needed contradiction.
\qed

\medskip

\section{The $(3,1)$-case}

The following provides a step towards establishing Theorem B:

\medskip

\noindent{\bf Proposition D} \, \it Let $\pi$ be an isobaric, algebraic representation of GL$_4({\mathbb A}_F)$ which is not of type $(3,1)$, with associated $4$-dimensional Hodge-Tate, $\overline{\mathbb Q}_p$-representation $\rho$ of ${\mathfrak G}_F$. Then $\rho^{\rm ss}$ is not of type $(3,1)$. Conversely, assuming in addition that $\pi$ is quasi-regular, if the decomposition type of $\rho$ is not $(3,1)$ for a sufficiently large $p$, then the isobaric type of $\pi$ is not $(3,1)$ either.
\rm

\medskip

\noindent{\it Proof}. \, Now suppose we have a decomposition as ${\mathfrak G}_F$-modules over $\overline {\mathbb Q}_p$:
\begin{equation}
\rho^{\rm ss} \, \simeq \, \tau \oplus \chi,
\label{dr:eqn3.1}
\end{equation}
 with $\tau$ (resp. $\chi$) of
dimension $3$ (resp. $1$). Since every subrepresentation of a
crystalline, resp. Hodge-Tate, representation is crystalline,
resp. Hodge-Tate, $\chi$ is crystalline. It follows that
$\chi$ is locally algebraic and by \cite{dr:Se}, it is defined by an
algebraic Hecke character $\chi$. Put
\begin{equation}
\nu \, = \, {\rm det}(\tau). \label{dr:eqn3.2}
\end{equation}
Then $\nu$ is also Hodge-Tate and corresponds to an algebraic
Hecke character~$\nu$.

Taking the contragredients of both sides of \eqref{dr:eqn3.1}, then twisting
by $\nu\chi^{-1}$, and noting that these processes
commute with taking semisimplification, we obtain
\begin{equation}
(\rho^\vee \otimes \nu\chi^{-1})^{\rm ss} \, \simeq
\, (\tau^\vee \otimes \nu\chi^{-1}) \oplus
\nu\chi^{-2}. \label{dr:eqn3.3}
\end{equation}

Appealing to hypothesis (b) of Theorem A, we see that outside a
finite set $S$ of places, the decompositions \eqref{dr:eqn3.1} and \eqref{dr:eqn3.3} imply
the following identity of $L$-functions:
\begin{equation}
L^S(s, \pi)L^S(s, \pi^\vee \otimes \nu\chi^{-1}) \, = \, L^S(s,
\tau)L^S(s, \tau^\vee \otimes
\nu\chi^{-1})L^S(s, \chi)L^S(s, \nu\chi^{-2}).
\label{dr:eqn3.4}
\end{equation}

\medskip

\noindent{\bf Lemma 3.1}\label{dr:lem3.1} \, \it Assuming \eqref{dr:eqn3.1}, we have
$$
\Lambda^2(\rho^{ss}) \, \simeq \, (\tau^\vee \otimes
\nu) \oplus (\tau \otimes \chi).
$$
\rm

\medskip

\noindent{\it Proof}. \, Using \eqref{dr:eqn3.1} and the fact that $\chi$
is one-dimensional and hence has trivial exterior square, we get
\begin{equation}
\Lambda^2(\rho^{\rm ss}) \, \simeq \, \Lambda^2(\tau)
\oplus \tau \otimes \chi. \label{dr:eqn3.6}
\end{equation}
So the Lemma will be proved if we verify that
\begin{equation}
\Lambda^2(\tau) \, \simeq \, \tau^\vee \otimes \nu.
\label{dr:eqn3.7}
\end{equation}
By the Tchebotarev density theorem, it suffices to check this easy
identity at the primes where all the representations are
unramified. Let $P$ such a prime, and let $\{\alpha_P,
\beta_P,\gamma_P\}$ denote the inverse roots of $\phi_P$, the
Frobenius at $P$, acting on $\tau$. Then the inverse roots on
$\Lambda^2(\tau)$ are given by the set
$$
\{\alpha_P\beta_P, \alpha_P\gamma_P, \beta_P\gamma_P\} \, = \,
\nu(\phi_P)\{\alpha_P^{-1}, \beta_P^{-1}, \gamma_P^{-1}\}.
$$
Here we have used the fact that $\nu$ is the determinant of
$\tau$. The identity \eqref{dr:eqn3.7} now follows because the inverse
roots of $\phi_P$ on $\tau^\vee$ are given by the inverses of
those on $\tau$. Done.

\qed

\medskip

This also finishes the proof of Proposition D.

\medskip

\noindent{\bf Proposition 3.2}\label{dr:prop3.2} \, \it Let $\rho$, $\pi$
satisfy the hypotheses of Theorem A. Then the decomposition \eqref{dr:eqn3.1}
cannot hold.
\rm

\medskip

\noindent{\it Proof}. \, Combining Lemma \ref{dr:lem3.1} with the
identity \eqref{dr:eqn3.4}, we get
\begin{equation}
L^S(s, \pi)L^S(s, \pi^\vee \otimes \nu\chi^{-1}) \, = \, L^S(s,
\Lambda^2(\rho^{\rm ss}) \otimes \chi^{-1})L^S(s,
\chi)L^S(s, \nu\chi^{-2}).\label{dr:eqn3.9}
\end{equation}

Now we appeal to a beautiful recent theorem of H.~Kim (\cite{dr:K}) which
establishes the automorphy in GL$(6)$ of the exterior square of
any cusp form on GL$(4)$. Applying this to our form $\pi$, we get
an isobaric automorphic representation $\Lambda^2(\pi)$ which is
functorial at all the unramified primes. Using this information in
\eqref{dr:eqn3.9}, and twisting by $\chi^{-1}$, we get
\begin{equation}
L^S(s, \pi \otimes \chi^{-1})L^S(s, \pi^\vee \otimes \nu\chi^{-2})
\, = \, L^S(s, \Lambda^2(\pi) \otimes \chi^{-1})\zeta^S(s)L^S(s,
\nu\chi^{-3}).\label{dr:eqn3.10} 
\end{equation}
(This identity slightly corrects (7.10) in \cite{dr:Ra1}, where $\chi^{-2}$ in
the second $L$-function on the left turns up mistakenly as $\chi^{-1}$, which however
does not affect the proof.)

One knows by Jacquet and Shalika (\cite{dr:JS1}) that for any isobaric
automorphic representation $\Pi_E$ of GL$(n)/{\mathbb Q}$ for any $n \geq 1$,
the incomplete $L$-function $L^S(s, \Pi_E)$ has no zero at $s=1$.
Consequently, $L^S(s, \Lambda^2(\pi) \otimes \chi^{-1})L^S(s,
\nu\chi^{-3})$ is non-vanishing at $s=1$. Since $\zeta^S(s)$ has a
pole at $s=1$, we see then that the right hand side, and hence the
left hand side, of \eqref{dr:eqn3.10} admits a pole at $s=1$. On the other
hand, since $\pi$ and $\pi^\vee$ are cusp forms on GL$(4)/{\mathbb Q}$, the
left hand side of \eqref{dr:eqn3.10} does not have a pole at $s=1$. So we get
a contradiction. The only possibility is that the decomposition of $\rho^{\rm ss}$
given by \eqref{dr:eqn3.1} cannot hold.

\qed

\medskip

Now let us prove the converse direction in Proposition~D. Suppose we have an isobaric decomposition
\begin{equation}
\pi \, \simeq \, \eta \boxplus \nu,\label{dr:eqn3.11}
\end{equation}
where $\eta$ is a cuspidal automorphic representation of GL$_3({\mathbb A}_F)$ and $\nu$ an idele class character of
$F$, necessarily with $\eta, \nu$ algebraic. Moreover, the quasi-regularity hypothesis on $\pi$
implies that $\eta$ is regular. Hence by the recent powerful theorem of Harris, Lan, Thorne and Taylor (\cite{dr:HLTT}), one may associate to $\eta$ a semisimple, $3$-dimensional $\overline {\mathbb Q}_p$-representation $\beta$. Since $\nu$ is algebraic, it also corresponds to an abelian $p$-adic representation $\nu$ of ${\mathfrak G}_F$. On the other hand, $\rho$ is associated to $\pi$. It follows by Tchebotarev that
\begin{equation}
\rho^{\rm ss} \, \simeq \, \beta \oplus \nu\label{dr:eqn3.12}
\end{equation}

Hence the Proposition will be proved if (the semisimple representation) $\beta$ is irreducible. Suppose not. If $\beta$ is a direct sum of three $p$-adic character, each of which is necessarily Hodge-Tate, we will deduce, by the strong multiplicity one theorem, that $\eta$ is isomorphic to an isobaric sum of three idele class characters of $F$, contradicting the cuspidality of $\eta$. So we may assume that we have
\begin{equation}
\beta \, \simeq \, \sigma \oplus \mu,\label{dr:eqn3.13}
\end{equation}
with $\sigma$ (resp. $\mu$) irreducible of dimension $2$ (resp. $1$). Let $\mu$ be the Hecke character attached to the (Hodge-Tate) $p$-adic character $\mu$.

Now the regularity of $\eta$ implies the same about the Hodge-Tate type of $\beta$. It follows that $\sigma$ ha distinct Hodge-Tate weights and by Taylor (\cite{dr:Ta2,dr:Ta3}), there is a finite solvable Galois extension $E$ of $F$ such that the restriction $\sigma_{ E}$ of $\sigma$ to ${\mathfrak G}_E$ is modular, i.e., is associated to a cusp form $\pi_0$ of GL$_2({\mathbb A}_E)$, which is regular algebraic. One can, by well known arguments via cyclic layers, descend $\pi_0$ to any subfield $M$ of $E$ with $E/M$ Galois and solvable and obtain a cusp form $\pi_0^M$ on GL$(2)/M$ which base changes to $\pi_0$ over $E$ and in fact corresponds to the restriction $\sigma_{M}$ of $\sigma$ to ${\mathfrak G}_M$.

As in \cite{dr:Ta2}, we now appeal to Brauer's theorem. By the inductive nature of $L$-functions, we have
\begin{equation}
L^{S^E}(s, \sigma_{E} \otimes \mu_{E}^{-1}) \, = \,
L^S(s, \sigma \otimes \mu^{-1} \otimes {\rm
Ind}_E^F(1_E)),
\label{dr:eqn3.14}
\end{equation}
where ${\rm Ind}_E^F(1_E))$ is the representation of
Gal$(\overline F/F)$ induced by the trivial representation $1_E$ of Gal$(\overline
F/E)$, $S$ is a finite set of places of $F$ containing the archimedean and ramified
places as well as the ones above $p$, and $S^E$ the finite set of places of $E$ above $S$. We can write
\begin{equation}
{\rm Ind}_E^F(1_E)) \, \simeq \, 1_F \, \oplus a_{E/F},
\label{dr:eqn3.15}
\end{equation}
for a unique representation $a_{E/F}$ of Gal$(\overline F/F)$
called the {\it augmentation representation}. Using Brauer's
theorem we can write it as a virtual sum of monomial
representations. More precisely,
\begin{equation}
a_{E/F} \, \simeq \, \oplus_{j=1}^r \,  n_j {\rm
Ind}_{E_j}^F(\alpha_j), \label{dr:eqn3.16}
\end{equation}
where for each $j$, $n_j$ is in ${\mathbb Z}$, $E_j$ a subfield of $E$ with
$E/E_j$ cyclic, and $\alpha_j$ a finite order character of
Gal$(\overline F/E_j)$ with trivial restriction to Gal$(\overline
F/E)$. Consequently, using \eqref{dr:eqn3.15}, \eqref{dr:eqn3.16}, the inductive nature and
additivity,
\begin{equation}
L^S(s, \sigma \otimes \mu^{-1} \otimes {\rm
Ind}_E^F(1_E)) \, = \, L^S(s, \sigma \otimes \mu^{-1})
\prod\limits_{j=1}^r \, L^{S_j}(s, \sigma_{E_j} \otimes
\mu_{E_j}^{-1} \otimes \alpha_j)^{n_j},
\label{dr:eqn3.17}
\end{equation}
where $S_j$ is the set of places of $E_j$ above $S$.
Now we can appeal to \eqref{dr:eqn3.14} and the modularity of $\sigma_{E_j}$ ($\forall \, j$)
to obtain
\begin{equation}
L^S(s, \sigma \otimes \mu^{-1}) \, = \, \prod_{j=0}^r \,
L^{S_j}(s, \pi_0^{E_j}\otimes \mu_{E_j}^{-1}\alpha_j)^{-n_j},
\label{dr:eqn3.18}
\end{equation}
where $n_0 = -1$, $E_0 = E$ and $\alpha_0=1$.

Since $\pi_0^{E_j}$ is cuspidal, $L^{S_j}(s, \pi_0^{E_j}\otimes \mu_{E_j}^{-1}\alpha_j)$ has no zero or pole at the right edge $s=1$. It then follows by \eqref{dr:eqn3.18} that
\begin{equation}
{\rm ord}_{s=1} \, L^S(s, \sigma \otimes \mu^{-1}) \, = \, 0.\label{dr:eqn3.19}
\end{equation}
Now applying \eqref{dr:eqn3.13} and the fact that $\beta$ is associated to $\eta$, we get
\begin{equation}
L^S(s, \eta\otimes\mu^{-1}) \, = \, L^S(s, \sigma\otimes\mu^{-1})\zeta_F^S(s).\label{dr:eqn3.20} 
\end{equation}
In view of \eqref{dr:eqn3.19}, the right hand side of \eqref{dr:eqn3.20} has a pole at $s=1$, which yields a contradiction as the left hand side (of \eqref{dr:eqn3.20}) is entire, $\eta$ being a cusp form on GL$(3)/F$.

This finishes the proof of Proposition D.
\qed

\medskip

\section{The $(2,2)$-case; End of Proof of Theorem B}

\medskip

We have just proved that $\pi$ is of type $(3,1)$ iff $\rho$ is also of the same type, under the running hypothesis that $\rho$ is crystalline and $p$ sufficiently large, i.e., with each of the Hodge-Tate weights of $\rho$ being $< 2(p-1)$. By a similar argument we can also prove that for $\pi$ Eisensteinian, $\rho^{\rm ss}$ has the same type as $\pi$. So we may {\it assume from here on that $\pi$ is cuspidal}, which is the key case. We already know that $\rho^{\rm ss}$ is not of type $(3,1)$, and it is easy to see that it cannot also be of type $(1,1,1,1)$ as then $\pi$ would, by the strong multiplicity one theorem, an isobaric sum of four Hecke characters.

Hence we will be done that it is impossible for $\rho^{\rm ss}$ to be of type $(2,2)$ or $(2,1,1)$ when $\pi$ is cuspidal. Suppose not. Then we have
\begin{equation}
\rho^{\rm ss} \, \simeq \, \sigma \oplus \sigma',\label{dr:eqn4.1}
\end{equation}
with dim$(\sigma)={\rm dim}(\sigma')=2$ and $\sigma$ irreducible. Thanks to Theorem A, both $\sigma$ and $\sigma'$ are odd. The Hodge-Tate types of $\sigma$ and $\sigma'$ are regular and they are crystalline as $\rho$ is by hypothesis. By Taylor (\cite{dr:Ta2,dr:Ta3}), we know that there is a finite Galois, totally real extension $E$ of $F$ such that $\sigma$ is modular over $E$, i.e., associated to a cusp form $\pi_0$ on GL$(2)/E$.

Suppose $\sigma'$ is reducible, i.e., of the form $\nu\oplus \nu'$ with dim$(\nu)={\rm dim}(\nu')=1$, then there are associated Hecke characters $\nu, \nu'$ of $F$, and by the argument at the end of the previous section, $L^S(s,\sigma\otimes\nu^{-1})$ is invertible at $s=1$. Consequently, the cuspidal $L$-function $L^S(s, \pi\otimes \nu^{-1})$, which equals $L^S(s,\sigma\otimes\nu^{-1})\zeta_F^S(s)L^S(s,\nu'/\nu)$, has a pole at $s=1$, which is a contradiction.

Hence we may assume that $\sigma, \sigma'$ are both irreducible and crystalline with $p$ sufficiently large. Since $\pi$ evidently satisfies the hypotheses of Theorem~A as well, we may apply Lemma~\ref{dr:lem1.3} and assume that $\pi$ remains cuspidal when base changed to any finite solvable normal extension which is totally real. Moreover, since the determinants $\nu, \nu'$ of $\sigma, \sigma'$ respectively are both odd by section 2, we may replace $F$ by an abelian, totally real extension over which $\nu=\nu'$, still with $\pi$ cuspidal.
Appealing to Taylor's theorem for $\sigma'$ as well, one gets a finite Galois, totally real extension $E'$ of $F$ and a cuspidal automorphic representation $\pi_0'$ of GL$_2({\mathbb A}_{E'})$ which is associated to $\sigma'_{E'}$. By Gelbart-Jacquet (\cite{dr:GJ}), there is an isobaric automorphic representation sym$^2(\pi_0)$, resp. sym$^2(\pi'_0)$, of GL$_3({\mathbb A}_E)$, resp. GL$_3({\mathbb A}_{E'})$ whose standard $L$-function equals $L^S(s, \pi_0; {\rm sym}^2)$, resp. $L^S(s, \pi'_0; {\rm sym}^2)$. Applying Brauer's theorem as in the previous section, and the modularity of $\sigma$, resp. $\sigma'$, over any sub $M \subset E$, resp. $M'\subset E'$, with $[E:M]$, resp. $[E':M']$, solvable and normal, we get the analogue of \eqref{dr:eqn3.19}:
\begin{equation}
{\rm ord}_{s=1} \, L^S(s, {\rm Ad}(\tau)) \, = \, 0, \, \, \, {\rm for} \, \, \, \tau\in\{\sigma, \sigma'\},\label{dr:eqn4.2}
\end{equation}
where ${\rm Ad}(\tau)={\rm sym}^2(\tau) \otimes \nu^{-1}$.
In addition, $L^S(s, {\rm Ad}(\tau))$ is meromorphic on the $s$-plane.

Applying the exterior square to \eqref{dr:eqn4.1} and twisting by $\nu^{-1}$, we obtain
\begin{equation}
L^S(s, \Lambda^2(\pi)\otimes \nu^{-1}) \, = \, L^S(\sigma\otimes\sigma'\otimes\nu^{-1})\zeta_F^S(s)^2,\label{dr:eqn4.3}
\end{equation}
which also gives the meromorphic continuation for $L^S(\sigma\otimes\sigma'\otimes\nu^{-1})$. Moreover, since $\pi$ is cuspidal, the left hand side of \eqref{dr:eqn4.3} can have at most a simple pole at $s=1$. So it follows that
\begin{equation}
{\rm ord}_{s=1} \, L^S(\sigma\otimes\sigma'\otimes\nu^{-1}) \, = \, 1.\label{dr:eqn4.4}
\end{equation}
On the other hand, applying the symmetric square to \eqref{dr:eqn4.1} and twisting, we get
\begin{equation}
L^S(s, \pi, {\rm sym}^2\otimes \nu^{-1}) \, = \, L^S(\sigma\otimes\sigma'\otimes\nu^{-1}) L^S(s, {\rm Ad}(\sigma))L^S(s, {\rm Ad}(\sigma')),\label{dr:eqn4.5}
\end{equation}
where the properties of the left hand side are given in \cite{dr:BuG} and \cite{dr:Sh1}. Moreover,
by \cite{dr:BuG}, $L^S(s, \pi, {\rm sym}^2\otimes \nu^{-1})$ has no zero at $s=1$. On the other hand, the right hand side of \eqref{dr:eqn4.5} must have a zero at $s=1$ by the conjunction of \eqref{dr:eqn4.2} and \eqref{dr:eqn4.4}. This gives the requisite contradiction to the decomposition \eqref{dr:eqn4.1}.

We have now proved that $\rho$ is irreducible when $\pi$ is cuspidal satisfying the hypotheses of Theorem B.

\qed

\section{Proof of Theorem C}

Let $\sigma, \sigma'$ be $2$-dimensional odd, semisimple, crystalline representations as in Theorem C, with $\sigma \oplus \sigma'$ associated to an isobaric automorphic form $\pi$ on GL$(4)/F$. If $\sigma, \sigma'$ are both reducible, they correspond to sums of algebraic Hecke characters, and the assertion of Theorem C follows by Hecke. More generally, if $\sigma, \sigma'$ become reducible over a finite Galois extension, then they are both of Artin-Hecke type, and the claim is well known (see \cite{dr:De}). So we may, and we will, assume that at least one of them, say $\sigma$, remains irreducible over any finite extension of $F$.

\medskip

\noindent{\bf Lemma 5.1}\label{dr:lem5.1} \, \it Under the hypothesis of irreducibility of $\sigma$, we must have
$$
\pi \, \simeq \, \eta \boxplus \eta',
$$
where $\eta, \eta'$ are isobaric automorphic forms on GL$(2)/F$, with $\eta$ cuspidal. Moreover, $\eta'$ is cuspidal iff $\sigma'$ is irreducible.
\rm

\medskip

Before proving this Lemma, let us note that if we knew $\pi$ to be algebraic and quasi-regular, we can conclude by the results in the earlier sections that $\pi$ cannot be of type $(3,1)$, $(2,1,1)$ or $(1,1,1,1)$. Similarly, if $\pi$ were regular and algebraic, we can also rule out $\pi$ being cuspidal.

\medskip

\noindent{\it Proof of Lemma 5.1} \, First suppose $\sigma'$ is reducible. Then evidently, since $\sigma$ is irreducible, $\sigma^\vee\otimes\sigma'$ has no Galois invariants. Moreover, applying the potential automorphy result of Taylor (\cite{dr:Ta1,dr:Ta2}) to the crystalline, odd $2$-dimensional $\sigma$, whose Hodge-Tate type satisfies his hypothesis relative to $p$, we deduce by Brauer (as it was done earlier) the following for any character $\mu$ of $F$:
\begin{equation}
{\rm ord}_{s=e}L^S(s, \sigma\otimes\mu) \, = \, 0,\label{dr:eqn5.2}
\end{equation}
where $e$ denotes the right edge (``Tate point''), which is $1+w/2$ if $\mu$ is of finite order. For suitable finite set $S$ of places of $F$, $L^S(s,\pi\otimes\mu)$ is a product of two abelian twists of the $L^S$-function of $\sigma$, which implies, by \eqref{dr:eqn5.2}, that $\pi$ must be an isobaric sum of the form $\eta \boxplus \mu\boxplus \mu'$, with $\eta$ a form on GL$(2)/F$ and $\mu, \mu'$ idele class characters of $F$. Moreover, $\eta$ must be cuspidal, for otherwise $\pi$ will be an isobaric sum of four characters, yielding an associated abelian $4$-dimensional representation of the Weil group of $F$, contradicting the irreducibility of $\sigma$. Putting $\eta'=\mu\boxplus\mu'$, we get, as asserted, $\pi \simeq \eta\boxplus \eta'$, with $\eta$ cuspidal and $\eta'$ not cuspidal.

We may now assume that $\sigma$ and $\sigma'$ are both irreducible, still with $\sigma$ remaining irreducible upon restriction to any open subgroup. Since they have the same parity and weight, we may consider the cyclic, totally real extension $E$ of $F$ such that their determinants are the same over $E$, say $\chi^w\nu$ with $\nu$ of finite order. Then the unitarily normalized $\pi_E$, the base change of $\pi$, has central character $\omega_E=\nu^2$ and satisfies $\pi_E^\vee \simeq \pi_E\otimes \nu^{-1}$. Comparing $L$-functions we get, for a finite set $T$ of places of $E$,
\begin{equation}
L^{T}(s, \Lambda^2(\pi_E)\otimes\nu^{-1}) \, = \, L^{T}(s, \sigma_E^\vee\otimes\sigma')\zeta_E^{T}(s)^2.\label{dr:eqn5.3}
\end{equation}
Suppose $\pi_E$ were cuspidal. Then the left hand side of \eqref{dr:eqn5.2} will have at most a simple pole at $s=1$, implying that $L^{T}(s, \sigma_E^\vee\otimes\sigma')$ must have a zero at $s=1$. On the other hand, we have
the identity \eqref{dr:eqn4.5} over $E$, which implies that $L^T(s, {\rm Ad}(\sigma_E))$ or $L^T(s, {\rm Ad}(\sigma'))$ must have a pole at $s=1$, since $L^T(s, \pi_E, {\rm sym}^2\otimes \nu^{-1})$ does not have a zero at $s=1$ by Shahidi (\cite{dr:Sh1,dr:Sh2}). Now we apply Taylor's theorem (\cite{dr:Ta1,dr:Ta2}) to conclude as before that $\sigma, \sigma'$ are potentially automorphic over finite Galois, totally real extensions $M, M'$ respectively of $F$, and we may deduce as before that for any character $\mu$ of $F$, and for any finite set $S$ of places of $F$,
\begin{equation}
{\rm ord}_{s=1}L^S(s, {\rm Ad}(\tau)\otimes\mu) \, = \, 0, \, \, \, {\rm for} \, \, \, \tau\in\{\sigma, \sigma'\}.\label{dr:eqn5.4}
\end{equation}
More precisely, we see that as $L^{T}(s, \sigma^\vee\otimes\sigma'\otimes\mu)$ has no zero at $s=1$,
\begin{equation}
-{\rm ord}_{s=1}L^{T}(s, \Lambda^2(\pi_E)\otimes\nu^{-1}) \, = \, -{\rm ord}_{s=1}L^{T}(s, \sigma_E^\vee\otimes\sigma'_E) + 2 \, \, \, \geq \, \, 2.\label{dr:eqn5.5}
\end{equation}
Similarly,
\begin{equation}
-{\rm ord}_{s=1}L^{T}(s, \pi_E; {\rm sym}^2\otimes\nu^{-1}) \, = \, -{\rm ord}_{s=1}L^{T}(s, \sigma_E^\vee\otimes\sigma'_E) \, \, \, \geq \, \, 0,\label{dr:eqn5.6}
\end{equation}
implying
\begin{equation}
-{\rm ord}_{s=1}L^{T}(s, \pi_E^\vee \times \pi_E) \, = \, -2{\rm ord}_{s=1}L^{T}(s, \sigma_E^\vee\otimes\sigma'_E) +2 \, \, \, \geq \, \, 2.\label{dr:eqn5.7}
\end{equation}
In particular, $\pi_E$ is not cuspidal. We claim that $\pi$ itself is not cuspidal. Suppose not. Then there exists a quadratic extension $M/L$ with $F\subset L\subset M\subset E$ such that $\pi_L$ is cuspidal, but $\pi_M$ is not, implying that $\pi_L$ admits a self twist by the quadratic character $\delta$, say, of $L$ associated to $M$. Then $\sigma_L\oplus \sigma'_L$ must also be isomorphic to its twist by $\delta$. Since $\sigma$ remains irreducible under restriction to any open subgroup, $\sigma$ remain irreducible, which forces an isomorphism $\sigma'_L \simeq \sigma_L\otimes \delta$. Then for non-trivial automorphism $\alpha$ of $L/F$, if any, $\sigma_L \simeq \sigma_L \otimes(\delta/\delta^\alpha)$, which implies that $\delta$ must be $\alpha$-invariant (since otherwise $\sigma_L$ would become reducible over $L(\delta/\delta^\alpha)$). It follows that already over $F$, $\sigma' \simeq \sigma\otimes\delta_0$, where $\delta_0$ is a descent of $\delta$ to $F$. Then $\sigma$ and $\sigma'$ have the same determinant, say $\nu\chi^w$, and (the assumption of cuspidality of $\pi$ implies)
$$
\sigma^\vee\otimes\sigma' \, \simeq \,  ({\rm Ad}(\sigma)\otimes\delta) \, \, \oplus \, \delta  \, \simeq \,
{\sigma'}^\vee\otimes\sigma.
$$
This implies by \eqref{dr:eqn5.4} that $L^S(s, \sigma^\vee\otimes\sigma')$ does not vanish at $s=1$. Consequently,
$L^S(\pi^\vee\times\pi)$, which equals $L^S(s, (\sigma\oplus\sigma')^\vee\otimes(\sigma\oplus\sigma'))$, has a pole of order $>1$ at $s=1$, which contradicts (by \cite{dr:JS1}) the cuspidality assumption on $\pi$. Thus $\pi$ is not cuspidal.
Moreover, if $\pi$ is of type $(3,1)$, then no abelian twist of its exterior square $L$-function can admit any pole at $s=1$. Thus $\pi$ must be of type $(2,2), (2,1,1)$ or $(1,1,1,1)$. As above we can eliminate $(1,1,1,1)$, and we can eliminate $(2,1,1)$ as well, for otherwise $L^S(s,\sigma\otimes\mu)L^S(s,\sigma'\otimes\mu)$ will have a pole at $s=e$ for a suitable character $\mu$, contradicting \eqref{dr:eqn5.2} (which holds for $\sigma'$ as well). Hence the assertion of the Lemma holds when $\sigma'$ is irreducible as well, with both $\eta, \eta'$ cuspidal.

\qed

\medskip

\noindent {\bf Proof of Theorem C (contd.)}  Comparing the exterior square $L$-functions of $\pi=\eta\boxplus\eta'$ and of $\sigma\oplus \sigma'$, we see that for any character $\mu$ of $F$,
\begin{equation}
L^S(s, \sigma^\vee\otimes\sigma'\otimes\mu) \, = \, L^S(s, \eta^\vee\boxtimes\eta'\otimes\mu),\label{dr:eqn5.8}
\end{equation}
where $\eta^\vee\boxtimes\eta'$ is the isobaric automorphic form on GL$(4)/F$ associated to the pair $(\eta^\vee, \eta')$ in \cite{dr:Ra3}.

Suppose $\sigma'$ is reducible. Then $\sigma\otimes\sigma'$ has no Galois invariants, since $\sigma$ is irreducible. On the other hand, as the right hand side of \eqref{dr:eqn5.8} is a product of two abelian twists of the standard $L$-function of the cusp form $\eta$, the order of $L^S(s, \sigma^\vee\otimes\sigma)$ at $s=1$ is $0$. So Theorem C is proved in this case.

It remains to consider when both $\sigma$ and $\sigma'$ are irreducible, in fact upon restriction to any open subgroup. Then $\sigma^\vee\otimes\sigma'$ has Galois invariants iff $\sigma\simeq \sigma'$, in which case the dimension of ${\mathfrak G}_F$-invariants is $1$. On the other hand, by Lemma~\ref{dr:lem5.1}, $\eta$ and $\eta'$ are both cuspidal on GL$(2)/F$. Then $L^S(s,\eta^\vee\boxtimes\eta')$ has a pole at $s=1$ iff $\eta$ is isomorphic to $\eta'$, in which case the pole is of order $1$. We get the same assertion for $L^S(s, \sigma^\vee\otimes\sigma'$ by \eqref{dr:eqn5.8}. So we have only to prove the {\it claim} that $\sigma\simeq \sigma'$ iff $\eta\simeq \eta'$. Arithmetically normalizing both sides, an isomorphism on either side furnishes an equality (at good primes $v$ of norm $q_v$) of the form:
$$
\{\alpha_v, q_v^w\alpha_v^{-1}, \beta_v, q_v^w\beta_v^{-1}\} \, = \, \{\gamma_v, q_v^w\gamma_v^{-1}, \gamma_v, q_v^w\gamma_v^{-1}\}.
$$
Up to renaming $\alpha_v$ as $\beta_v$ and switching $\gamma_v$ with $q_v^w\gamma_v^{-1}$, we may assume that $\alpha_v=\gamma_v$. Then $q_v^w\alpha_v^{-1}$ equals $q_v^w\gamma_v^{-1}$, resulting in the equality
$\{\beta_v, q_v^w\beta_v^{-1}\} =\{\gamma_v, q_v^w\gamma_v^{-1}\}$. The claim follows.

\qed



\end{document}